\documentclass{amsart}
\usepackage{amsmath,enumerate}

\newtheorem{thm}{Theorem}

\newtheorem{cor}[thm]{Corollary}
\newtheorem{lem}[thm]{Lemma}
\newtheorem{defi}[thm]{Definition}
\newtheorem{remark}[thm]{Remark}
\newtheorem{example}[thm]{Example}
\newtheorem{pb}[thm]{Problem}

\newcommand{\real}{{\mathbb R}}
\newcommand{\nat}{{\mathbb N}}

\newcommand{\comp}{{\mathbb C}}

\newcommand{\E}{{\mathcal E}}

\newcommand{\M}{{\mathcal M}}

\renewcommand{\a}{\alpha}

\newcommand{\e}{\varepsilon}
\renewcommand{\th}{\theta}

\renewcommand{\e}{\varepsilon}
\newcommand{\f}{\varphi}

\newcommand{\s}{\sigma}

\newcommand{\tr}{\mbox{\rm tr}}
\newcommand{\ot}{\otimes}

\renewcommand{\t}{\tau}
\newcommand{\8}{\infty}
\newcommand{\el}{\ell}

\newcommand{\wt}{\widetilde}
\newcommand{\wh}{\widehat}
\newcommand{\n}{\noindent}
\newcommand{\pf}{\noindent{\it Proof.~~}}
\newcommand{\cqd}{\hfill$\Box$}
\newcommand{\be}{\begin{eqnarray*}}
\newcommand{\ee}{\end{eqnarray*}}
\newcommand{\beq}{\begin{equation}}
\newcommand{\eeq}{\end{equation}}

\title{A description of $\big(C_p[L_p(M)],\; R_p[L_p(M)]\big)_\th$ }
\author{Quanhua Xu}
\address{D\'{e}partement de Math\'{e}matiques\\
Universit\'{e} de France-Comt\'{e}\\
16 Route de Gray\\
25030 Besan\c con Cedex,  France } \email{qx@math.univ-fcomte.fr}

\begin{document}

 \maketitle


\begin{abstract}
We give a simple explicit description of the norm in the complex
interpolation space $(C_p[L_p(M)],\; R_p[L_p(M)])_\th$ for any von
Neumann algebra $M$ and any $1\le p\le\8$.
\end{abstract}

\vskip 1.5cm


 Let $M$ be a semifinite von Neumann algebra equipped with a
normal semifinite faithful trace $\t$. Let $L_p(M)$ be the
associated noncommutative $L_p$-space. Given an integer $n$ let
$C^n[M]$ (resp. $R^n[M]$) be $M^n$ equipped with the following
norm
 $$\big\|\sum_{k=1}^n x_k^*x_k\big\|_M^{1/2}
 \quad \big(\mbox{resp.}\quad
 \big\|\sum_{k=1}^n x_kx_k^*\big\|_M^{1/2}\,\big)\;.$$
We view $(C^n[M],\;R^n[M])$ as a compatible couple by identifying
algebraically both $C^n[M]$ and $R^n[M]$ with $M^n$. Then we can
consider the complex interpolation space $(C^n[M],\; R^n[M])_\th$
for $0<\th<1$ (cf. \cite{bl} for complex interpolation). Pisier
\cite{pis-proj} described this interpolation norm by the following
simple formula: For any $(x_1, ..., x_n)\in M^n$
 \beq\label{pisier}
 \big\|(x_1, ..., x_n)\big\|_{(C^n[M],\; R^n[M])_\th}
 =\big\|\sum_{k=1}^n L_{x_k^*}R_{x_k}\big\|_{B(L_p(M))}\;,
 \eeq
where $1/p=\th$, and where $L_x$ (resp. $R_x$) denotes the
multiplication on $L_p(M)$ by $x$ from the left (resp. right).
Haagerup \cite{haag-self} then extended this formula to any von
Neumann algebra $M$, at least for $\th=1/2$. In this case $L_2(M)$
can be any Hilbert space at which $M$ acts standardly.

Pisier used (\ref{pisier}) as a tool in his study of the problem
when there is a contractive projection from a super von Neumann
algebra $N$ onto $M$. In particular, when $N=B(L_2(M))$ ($M$ being
then represented by left multiplication on $L_2(M)$), this problem
reduces to the injectivity of $M$.  We refer the interested reader
to \cite{pis-oh}, \cite{pis-proj}, \cite{pis-edin} and
\cite{haag-self} for more information. Here we content ourselves
only by mentioning the following result from \cite{pis-oh}: $M$ is
injective iff
 $$\big\|(x_1, ..., x_n)\big\|_{(C^n[M],\; R^n[M])_{1/2}}
 =\big\|\sum_{k=1}^n x_k\ot \bar x_k
 \big\|_{M\bar\ot_{\min}\overline M}\;.$$

The purpose of this note is to consider the $L_p$-space version of
(\ref{pisier}) for any $1\le p\le\8$, i.e. to give a simple
description of the norm of the interpolation space
$(C^n[L_p(M)],\; R^n[L_p(M)])_\th$ for any von Neumann algebra
$M$.

 It is well known by now that there are several equivalent
constructions of the noncommutative $L_p$-spaces associated with a
general von Neumann algebra. In this note we use  those
constructed by Haagerup \cite{haag-Lp}. Our reference for these
spaces is \cite{terp-Lp}. In the sequel, $M$ will denote a general
von Neumann algebra, unless explicitly stated otherwise. $L_p(M)$
stands for the Haagerup noncommutative $L_p$-space based on $M$.
However, whenever $M$ is semifinite, we will always consider
$L_p(M)$ as defined  by a normal semifinite faithful trace. We
refer to the survey \cite{px-survey} for semifinite noncommutative
$L_p$-spaces and references therein.

Let $M$ be a von Neumann algebra and $1\le p\le\8$. Given an
integer $n\in\nat$ we denote by $C^n_p[L_p(M)]$ (resp.
$R^n_p[L_p(M)]$) $L_p(M)^n$ equipped with the following norm
 $$\big\|\big(\sum_{k=1}^n x_k^*x_k\big)^{1/2}\big\|_p
 \quad \big(\mbox{resp.}\
 \big\|\big(\sum_{k=1}^n x_kx_k^*\big)^{1/2}\big\|_p\,\big)\;.$$
$C^n_\8[L_\8(M)]$ (resp. $R_\8^n[L_\8(M)]$) is, of course,
$C^n[M]$ (resp. $R^n[M]$) introduced previously. As before in the
case of $p=\8$ we regard $(C^n_p[L_p(M)],\; R^n_p[L_p(M)])$ as a
compatible couple by identifying $C^n_p[L_p(M)]$ and
$R^n_p[L_p(M)]$ with $L_p(M)^n$. The main result of this note is
the following generalization of (\ref{pisier}). $\|\;\|_p$ denotes
the norm in $L_p(M)$.

\begin{thm}\label{interp}
 Let $1\le p\le\8$ and $0<\th<1$. Let $r, r_0(\th)$ and $r_1(\th)$ be
determined by
 $$\frac{1}{r}=1-\frac{2}{\max(p, p')}\;,\quad
 \frac{1}{r_0(\th)}=\frac{\th}{2r}\;,\quad
 \frac{1}{r_1(\th)}=\frac{1-\th}{2r}\;,$$
where $p'$ is the index conjugate to $p$. Let $x=(x_1, ...,
x_n)\in L_p(M)^n$.
\begin{enumerate}[{\rm i)}]
 \item If $p\le2$ , then
 \beq\label{interp1}
 \|x\|_{(C^n_p[L_p(M)],\; R^n_p[L_p(M)])_\th}=
 \inf\Big\{\|a\|_{r_0(\th)}\,
 \|b\|_{r_1(\th)}\,\big(\sum_k\|y_k\|^2_2\big)^{1/2}\Big\},
 \eeq
where the infimum runs over all factorizations of $x$ as
$x_k=ay_kb$ with $a\in L_{r_0(\th)}(M)$, $b\in L_{r_1(\th)}(M)$)
and $y_k\in L_2(M)$ $(1\le k\le n)$.
 \item If $p\ge 2$, then
 \beq\label{interp2}
 \|x\|_{(C^n_p[L_p(M)],\; R^n_p[L_p(M)])_\th}=
 \sup\Big\{\big(\sum_k\|a x_k b\|_2^2\big)^{1/2}\Big\},
 \eeq
where the supremum runs over all $a$ and $b$ respectively in the
unit balls of $L_{r_0(\th)}(M)$ and $L_{r_1(\th)}(M)$.
\end{enumerate}
\end{thm}

After having completed this note,  we learnt from Marius Junge
that he and Parcet had obtained a result similar to (even more
general than) Theorem \ref{interp} (see \cite{jp}).

\medskip

Clearly, (\ref{interp2}) in the case of $p=\8$ reduces to
(\ref{pisier}) for a semifinite $M$. For a general $M$ we get the
following extension of Haagerup's result to all $\th\in(0,1)$.

\begin{cor}
 Let $M$ be a von Neumann algebra. Let $0<\th<1$ and
$\frac{1}{p}=\th$. Then for any $x=(x_1,..., x_n)\in M^n$ we have
  \be
 \big\|(x_1, ..., x_n)\big\|_{(C^n[M],\; R^n[M])_\th}
 =\big\|\sum_{k=1}^n L_{x_k^*}R_{x_k}\big\|_{B(L_p(M))}\;.
 \ee
\end{cor}

It is a routine exercice to extend the theorem above to the case
of infinite sequences. Indeed, let $C_p[L_p(M)]$ be the completion
(relative to the w*-topology for $p=\8$) of the family of all
finite sequences in $L_p(M)$ with respect to the following norm
 $$\big\|\big(\sum_n x_n^*x_n\big)^{1/2}\big\|_p\;.$$
It is easy to see that  $C_p[L_p(M)]$ consists of all sequences
$x=(x_n)$ in $L_p(M)$ such that
 $$\sup_n\,\big\|\big(\sum_{k=1}^n x_k^*x_k
 \big)^{1/2}\big\|_p<\8$$
and the norm of $x$ is equal to the supremum above. Similarly, we
define $R_p[L_p(M)]$ as the space of all sequences $(x_n)$ in
$L_p(M)$ such that $(x_n^*)\in C_p[L_p(M)]$ equipped with the norm
 $$\|(x_n)\|_{R_p[L_p(M)]}=\|(x^*_n)\|_{C_p[L_p(M)]}\;.$$
It should be pointed out that the two norms in $C_p[L_p(M)]$ and
$R_p[L_p(M)]$ are in general not comparable at all. Again, we
view $(C_p[L_p(M)],\;R_p[L_p(M)])$ as a compatible couple by
injecting both $C_p[L_p(M)]$ and $R_p[L_p(M)]$ into
$\el_\8(L_p(M))$. Then Theorem \ref{interp} still holds for
$(C_p[L_p(M)],\;R_p[L_p(M)])_\th$ without any change, except that
in the case of $p=\8$, the norm on the left hand side of
(\ref{interp2}) should be replaced by that of
$(C_\8[L_\8(M)],\;R_\8[L_\8(M)])^\th$, the space constructed by
the second complex interpolation method

\medskip

Before proceeding to the proof of Theorem \ref{interp}, let us
make some comments for the readers familiar with operator space
theory. First, such a reader might have already realized that
$C_p[L_p(M)]$ (resp. $R_p[L_p(M)]$) is not only a pure notation
but the $p$-column space $C_p$ (resp. the $p$-row space $R_p$)
with values in $L_p(M)$ in Pisier's language \cite{pis-ast}. Here,
$C_p$ (resp. $R_p$) is defined as the (first) column (resp. row)
subspace of the Schatten class $S_p$. All noncommutative
$L_p$-spaces are equipped with their natural operator space
structure (see \cite{pis-ast}, \cite{pis-intro} and
\cite{ju-fubini}).

Second, let $M$ be an injective von Neumann algebra and $E$ an
operator space. We then have the vector-valued noncommutative
$L_p$-space $L_p[M;E]$ as defined in \cite{pis-ast}. (Note that
this is done in \cite{pis-ast} with the additional assumption that
$M$ is semifinite. However, the type III case can be dealt with
similarly.) In this language we have
 $$C_p[L_p(M)]=L_p[M; C_p]\quad\mbox{and}\quad
 R_p[L_p(M)]=L_p[M; R_p].$$
Then by \cite{pis-ast}, for any $0<\th<1$
 $$(L_p[M; C_p],\;L_p[M; R_p]_\th=L_p[M; (C_p,\; R_p)_\th]
 =L_p[M; C_q],$$
where $\frac{1}{q}=\frac{1-\th}{p}+\frac{\th}{p'}$. Thus for an
injective $M$, Theorem \ref{interp} can  be restated as follows.

\begin{cor}\label{interpbis}
 Keep the notations in Theorem \ref{interp} with
the additional assumption that $M$ is injective. Let $q$ be
defined by $\frac{1}{q}=\frac{1-\th}{p}+\frac{\th}{p'}$.  Let $x$
be a finite sequence in $L_p(M)$. Then if $p\le 2$,
 \beq\label{interp1bis}
 \|x\|_{L_p[M; C_q]}=
 \inf\Big\{\|a\|_{r_0(\th)}\,
 \|b\|_{r_1(\th)}\, \big(\sum_k\|y_k\|^2_2\big)^{1/2}\Big\},
 \eeq
where the infimum runs over all factorizations of $x$ as
$x_k=ay_kb$  with $a\in L_{r_0(\th)}(M)$, $b\in L_{r_1(\th)}(M)$
and $y_k\in L_2(M)$ $(k\ge1)$; if $p\ge 2$,
 \beq\label{interp2bis}
 \|x\|_{L_p[M; C_q]}=
 \sup\Big\{\big(\sum_k\|a x_k b\|_2^2\big)^{1/2}\Big\},
 \eeq
where the supremum runs over all $a$ and $b$ respectively in the
unit balls of $L_{r_0(\th)}(M)$ and $L_{r_1(\th)}(M)$.
\end{cor}

\n{\bf Remark.} In the case where $M=B(\el_2)$, (\ref{interp2bis})
has been already known before. In fact, the case $p=\8$ is
\cite[Theorem 8.4]{pis-oh}. On the other hand, the case $2\le
p<\8$ is proved in \cite{xu-gro}. These results are particularly
useful for studying the complete boundedness of maps with values
in $C_q$ (see \cite{pis-Rp}, \cite{xu-gro}).

\medskip

 The rest of the note is devoted to the proof of
Theorem \ref{interp}. For notational simplicity, we will denote
the norm of $(C^n_p[L_p(M)],\; R^n_p[L_p(M)])_\th$ by $\|\;\|_{p,
\th}$, and given $x=(x_1,..., x_n)\in (L_p(M))^n$ define $\a_{p,
\th}(x)$ to be the infimum in (\ref{interp1}) if $p\le 2$ and the
supremum in (\ref{interp2}) if $p\ge 2$. In the following, we will
also use these notations for $\th=0$ and $\th=1$. In these extreme
cases, $\a_{p, \th}(x)$ is defined by the same formulas; but
$\|x\|_{p,0}$ (resp. $\|x\|_{p,1}$) must be, of course, replaced
by $\|x\|_{C^n_p[L_p(M)]}$ (resp. $\|x\|_{R^n_p[L_p(M)]}$), as
usual in interpolation theory.

We will need the following lemma.

\begin{lem}\label{lemma}
 Let $2<p<\8$ and $p'$ be the conjugate index of $p$. Then the
dual norm of $\a_{p, \th}$ on $L_p(M)^n$ is equal to $\a_{p',
\th}$.
\end{lem}

\pf Let $y=(y_1,..., y_n)\in L_{p'}(M)^n$. Define
 $$\el_y\;:\; L_p(M)^n\to\comp\quad\mbox{by}\quad
 \el_y(x)=\sum_k\tr(y_k^*x_k)\;,$$
where $\tr$ is the distinguished tracial functional on $L_1(M)$.
Then it is easy to see that
 $$|\el_y(x)|\le\a_{p', \th}(y)\,\a_{p, \th}(x).$$
It follows that $\el_y$ is a continuous functional on
$(L_p(M)^n,\; \a_{p, \th})$ and of norm $\le\a_{p', \th}(y)$.

Conversely, assume that $\el$ is a continuous functional on
$(L_p(M)^n,\; \a_{p, \th})$. Since the restriction of $\a_{p,
\th}$ to each component of $L_p(M)^n$ coincides with the norm of
$L_p(M)$, there is $y\in L_{p'}(M)^n$ such that $\el=\el_y$ as
above. Thus for any $x\in L_{p}(M)^n$ we have
 \beq\label{lem1}
 |\el(x)|=\big|\sum_k\tr(y_k^*x_k)\big|\le \|\el\|\,\a_{p, \th}(x)
 =\|\el\|\,\sup\big\{\big|\sum_k\tr(z_k^*ax_kb)\big|\big\},
 \eeq
where the supremum is taken over all $a, b$ and $z$, respectively,
in the unit balls of $L_{r_0(\th)}(M)$, $L_{r_1(\th)}(M)$ and
$\el_2^n(L_2(M))$. One can further require $a$ and $b$ to be
positive. Since the left and right supports of the $y_k$ are
$\s$-finite projections, so is their supremum $e$. Replacing $M$
by the reduced algebra $eMe$ if necessary, we can assume $M$
itself  $\s$-finite, and so $L_p(M)$ can be constructed from a
normal faithful state on $M$.

By a typical minimax principle (cf. e.g. \cite[Lemma
2.3.1]{er-book}), or alternatively, a Hahn-Banach separation
argument by convexifying (\ref{lem1}) (cf. \cite{pis-gro}) , we
deduce from (\ref{lem1}) that there are $a\ge 0, b\ge 0$ and $z$ ,
respectively, in the unit balls of $L_{r_0(\th)}(M)$,
$L_{r_1(\th)}(M)$ and $\el_2^n(L_2(M))$ such that
 $$
 |\el(x)|\le\|\el\|\,\big|\sum_k\tr(z_k^*ax_kb)\big|
 \le\|\el\|\,\big(\sum_k\|ax_kb\|_2^2\big)^{1/2},
 \quad \forall\; x\in L_p(M)^n\;.
 $$
Hence $(ax_1b,..., ax_nb)\mapsto \el(x)$ extends to a continuous
functional on $\el_2^n(L_2(M))$ of norm $\le\|\el\|$. Therefore,
there is $u\in \el_2^n(L_2(M))$ such that
 $$\el(x)=\sum_k\tr(u_k^*ax_kb)\quad
 \mbox{and}\quad\big(\sum_k\|u_k\|_2^2\big)^{1/2}\le\|\el\|.$$
Recalling that $\el=\el_y$, we have
 $$\tr(y_k^*x_k)=\tr(u_k^*ax_kb)=\tr(bu_k^*ax_k),\quad
 \forall\;x_k\in L_p(M),\; 1\le k\le n.$$
It then follows that $y_k=au_kb$ for all $k$, and so
$\a_{p',\th}(y)\le\|\el\|$. Thus the lemma is proved \cqd

\medskip

\n{\bf Remarks.}  i) The proof above also applies to normal
functionals on $(M^n, \a_{\8, \th})$: if $\el$ is a normal
functional on $(M^n, \a_{\8, \th})$, then $\el=\el_y$ for some
$y\in L_1(M)^n$ and $\|\el\|=\a_{1,\th}(y)$. Consequently,
$(L_1(M)^n,\; \a_{1,\th})$ is the predual of $(M^n, \a_{\8,
\th})$.

ii) Let $1\le q<2$. It is easier to show that the dual space of
$(L_q(M)^n,\; \a_{q,\th})$ is $(L_{q'}(M)^n,\; \a_{q',\th})$. Then
we can recover Lemma \ref{lemma} by reflexivity. However, this
argument does not seem to yield the previous remark on the normal
functionals on $(M^n, \a_{\8, \th})$, which will be also needed
later.

\medskip

We will further need a well-known unpublished result by Haagerup
\cite{haag-red}. To state it, we first recall that a weight $\f$
on $M$ is called strictly normal if there is a family
$\{\f_i\}_{i\in I}$ of normal positive functionals with pairwise
disjoint supports such that
 $$\f=\sum_i\f_i\;.$$
Any von Neumann algebra admits a strictly normal semifinite
faithful weight.  As usual,  $\s_t^\f$ stands for the modular
automorphism group of a weight $\f$.

\begin{thm}\label{haag} {\bf (Haagerup)}
 Let $\f$ be a strictly normal semifinite faithful weight on
$M$. Then there are a von Neumann algebra $\M$, a strictly normal
semifinite faithful weight $\wh\f$ on $\M$ and an increasing
family $\{ \M_i \}_{i\in I}$ of w*-closed $\ast$-subalgebras of
$\M$ satisfying the following properties :
 \begin{enumerate}[{\rm i)}]
 \item $M$ is a von Neumann subalgebra of $\M$,
$\wh\f\big|_M=\f$ and $\s_t^{\wh \f}\big|_M=\s_t^\f$ for all
$t\in\real$;
 \item there is a normal faithful conditional expectation $\E$
 from $\M$ onto $M$ such that $\wh\f\circ\E=\f$ and
 $\s_t^{\wh\f}\circ\E=\E\circ\s_t^{\wh\f}$ for all $t\in\real$;
 \item  each $\M_i$ is finite and $\sigma$-finite and their union
is w*-dense in $\M$;
 \item  for every $i\in I$ there is a normal conditional
expectation $\E_i$ from $\M$ onto $\M_i$ such that
$$\E_i\circ\E_j=\E_j\circ\E_i\ \mbox{\rm whenever}\ i\le j\quad
\mbox{\rm and}\quad \s_t^{\wh\f}\circ\E_i= \E_i\circ
\s_t^{\wh\f},\quad t\in\real,\ i\in I.$$
 \end{enumerate}
 \end{thm}

Now we are in a position to prove Theorem \ref{interp}.

\medskip

 \n{\em Proof of Theorem \ref{interp}.} Both
(\ref{interp1}) and (\ref{interp2}) are trivially true for $p=2$.
Thus in the following we assume $p\not=2$.  The proof is divided
into three steps. We will first prove (\ref{interp1}) in the case
where $M$ is semifinite (or only finite). The main ingredient of
this part is an operator-valued version of the classical Szeg\"o
factorization theorem.  Then we will use Theorem \ref{haag} to
reduce the general case to the finite one. Finally, we will get
(\ref{interp2}) by duality from (\ref{interp1}).

\medskip

\n{\em Step~1: The proof of (\ref{interp1}) for finite $M$.} Let
$1\le p< 2$. In this first part we will prove (\ref{interp1}) with
the additional assumption that $M$ is finite and $\sigma$-finite.
Thus we can assume that $L_p(M)$ is defined on $M$ by a normal
finite faithful normalized trace $\t$.

Let us first prove that $\|x\|_{p,\th}\le \a_{p,\th}(x)$ for any
$x=(x_1,..., x_n)\in (L_p(M))^n$. This is easy for $\th=0$ and
$\th=1$. Indeed, let $ \a_{p,0}(x)<1$, and let $x_k=ay_kb$ be a
factorization of $x$ such that
 $$\|a\|_\8<1,\quad  \|b\|_{2r}<1\quad  \mbox{and}\quad
 \sum_{k=1}^n\|y_k\|_2^2<1\;.$$
Then
 $$\sum_k x_k^*x_k\le \|a\|_\8^2\,\sum_kb^*\,y_k^*y_k\,b
 \le b^*\,\big(\sum_k y_k^*y_k\big)b\,.$$
Thus by the H\"older inequality
 $$\|x\|_{p,0}=\big\|\big(\sum_k x_k^*x_k\big)^{1/2}\big\|_p
 \le \|b^*\|_{2r}^{1/2}\,\|b\|_{2r}^{1/2}\, \big\|\big(\sum_k
 y_k^*y_k\big)^{1/2}\big\|_2<1.$$
It follows that $\|x\|_{p,0}\le \a_{p,0}(x)$. Similarly,
$\|x\|_{p,1}\le \a_{p,1}(x)$. Then by complex interpolation for
trilinear maps (cf. \cite[Theorem 4.1.1]{bl}), we deduce
$\|x\|_{p,\th}\le \a_{p,\th}(x)$ for all $0<\th<1$.

\medskip

It is the converse inequality which is non trivial. The following
proof is similar to the proof of \cite[Theorem~2.1]{pis-proj}. Fix
an $x\in (L_p(M))^n$ such that $\|x\|_{p,\th}<1$. Let
$S=\{z\in\comp\;:\; 0\le{\rm Re}(z)\le1\}$. Let $\partial_0$ and
$\partial_1$ be respectively the right and left border of $S$.
Then there is a continuous function $F: S\to (L_p(M))^n$ such that
$F$ is analytic in the interior of $S$, $F(\th)=x$, and such that
 $$\sup_{z\in\partial_0} \|F(z)\|_{C^n_p[L_p(M)]}<1,\quad
 \sup_{z\in\partial_1} \|F(z)\|_{R^n_p[L_p(M)]}<1.$$
Write $F=(F_1, ..., F_n)$ and let $\e$ be a fixed (small) positive
number. Define
 $$X(z)=\big(\e + \sum_k F_k(z)^*F_k(z)\big)^{\frac{1}{2}}\ \mbox{for}\
 z\in\partial_0\quad\mbox{and}\quad
 X(z)=\big(\e + \sum_k F_k(z)F_k(z)^*\big)^{\frac{1}{2}}
 \ \mbox{for}\ z\in\partial_1\;.$$
Note that $X(z)$ is a positive invertible measurable operator with
$X(z)^{-1}\in M$ for every $z\in\partial_0\cup \partial_1$. Since
 $$F_k(z)^*\, F_k(z)\le X(z)^2,\quad z\in\partial_0,$$
there is $u_k\;:\; \partial_0\to M$ such that
 \beq\label{pf1}
 F_k(z)=u_k(z)X(z)\quad\mbox{and}\quad
 \sum_k u_k(z)^* u_k(z)\le 1,\quad z\in\partial_0.
 \eeq
Similarly, there is $v_k\;:\; \partial_1\to M$ such that
 $$F_k(z)=X(z)v_k(z)\quad\mbox{and}\quad
 \sum_k v_k(z)v_k(z)^*\le 1,\quad z\in\partial_1.$$
Define
 \beq\label{pf2}
 A(z)= 1\ \mbox{for}\ z\in\partial_0
 \quad\mbox{and}\quad
 A(z)=X(z)^{1-\frac{p}{2}} \ \mbox{for}\ z\in\partial_1;
 \eeq
 $$B(z)=X(z)^{1-\frac{p}{2}}  \ \mbox{for}\ z\in\partial_0
 \quad\mbox{and}\quad
 B(z)= 1 \ \mbox{for}\ z\in\partial_1;$$
 \beq\label{pf3}
 W_k(z)=u_k(z)X(z)^{\frac{p}{2}}\ \mbox{for}\ z\in\partial_0
 \quad\mbox{and}\quad
 W_k(z)=X(z)^{\frac{p}{2}}\,v_k(z) \ \mbox{for}\ z\in\partial_1.
 \eeq
Then we have the following factorization
 \beq\label{pf4}
 F_k(z)=A(z)W_k(z)B(z), \quad z\in \partial_0\cup
 \partial_1.
 \eeq

Now we use a well-known  Szeg\"o type factorization for
operator-valued analytic functions to bring the factorization
(\ref{pf4}) to an analytic one. The result we need here is
\cite[Corollary 8.2]{px-survey} applied to the special case of
Example (iii) in \cite[p.1496]{px-survey}.  We should emphasize
that this result is a combination (as well as a certain
improvement) of several previous results due to notably
Helson-Lowdenslager, Winner-Masani, Devinatz and Sarason. We refer
to \cite{px-survey} for more information and more historic
references. We also note that for our purpose we can use instead
\cite{xu-studia} plus an approximation argument.

Thus by \cite[Corollary 8.2]{px-survey}, there are two functions
$\Phi$ and $\Psi$ defined on $S$ with values in $L_{2r}(M)$ such
that $\Phi$ and $\Psi$ are analytic in the interior of $S$, and
such that
 \beq\label{pf5}
 \Phi(z)\Phi(z)^*=A(z)^2\quad\mbox{and}\quad
 \Psi(z)^*\Psi(z) =B(z)^2,\quad z\in\partial_0\cup
 \partial_1.
 \eeq
Moreover, both $\Phi(z)$ and $\Psi(z)$ are invertible with bounded
inverses in $M$ for all $z\in S$. Instead of \cite[Corollary
8.2]{px-survey}, we can directly use \cite[Theorem
8.12]{px-survey} (which is due to Saito). Indeed, let
$w(z)=B(z)^{-1}$. Then $w(z)\in M$ and $w(z)^{-1}\in L_{2r}(M)$
for every $z\in \partial_0\cup
 \partial_1$. Since $r\ge 1$ and $M$ is finite, $L_{2r}(M)\subset
 L_{2}(M)$. Thus by \cite[Theorem
8.12]{px-survey}, there are a function $\wt u$ such that $\wt
u(z)\in M$ is unitary  and an invertible analytic function $\psi$
such that $w(z)=\wt u(z)\psi(z)$ for every $z\in \partial_0\cup
 \partial_1$. Set $\Psi (z)=\psi(z)^{-1}$. Then $\Psi$ is    an
 invertible analytic function and $B(z)=\wt u(z)\Psi(z)$. It then
 follows that $\Psi(z)^*\Psi(z) =B(z)^2$ for all $z\in\partial_0\cup
 \partial_1$, as required.

(\ref{pf5}) implies that there are $U\;:\; S\to M$ and $V\;:\;
S\to M$ such that
 \beq\label{pf6}
 A(z)=\Phi(z)U(z),\quad B(z)=V(z)\Psi(z),\quad
 \|U(z)\|_\8\le 1,\quad \|V(z)\|_\8\le 1,\quad
 z\in\partial_0\cup\partial_1.
 \eeq
Thus
 \be
 F_k(z)=\Phi(z)\big[U(z)W_k(z)V(z)\big]\Psi(z)
 \buildrel{\rm def}\over=\Phi(z)Y_k(z)\Psi(z).
 \ee
Since $\Phi(z)^{-1}$ and $\Psi(z)^{-1}$ are analytic in the
interior of $S$, so is $Y_k(z)=\Phi(z)^{-1}F_k(z)\Psi(z)^{-1}$ for
every $k$. Therefore, we have the desired analytic factorization.

Let us estimate the norms of each factor on the border of $S$. By
the choice of $\Phi(z)$ in (\ref{pf5}) and the definition of
$A(z)$ in (\ref{pf2}), we have
 $$\sup_{z\in\partial_0}\|\Phi(z)\|_\8=
 \sup_{z\in\partial_0}\|A(z)\|_\8=1$$
and
 $$\sup_{z\in\partial_1}\|\Phi(z)\|^{2r}_{2r}=
 \sup_{z\in\partial_1}\|A(z)\|^{2r}_{2r}=
 \sup_{z\in\partial_1}\|X(z)\|^p_p<1,$$
provided $\e$ is small enough. Therefore, by interpolation
 $$\|\Phi(\th)\|_{r_0(\th)}\le 1.$$
Similarly, $\|\Psi(\th)\|_{r_1(\th)}\le 1$.

Concerning $Y_k$,  for any $z\in\partial_0$ by (\ref{pf6}), the
definition of $W_k$ in (\ref{pf3}) and the inequality in
(\ref{pf1}), we have
 $$\sum_k\|Y_k(z)\|_2^2\le \sum_k\|W_k(z)\|_2^2
 =\big\|X(z)^{\frac{p}{2}}\,\sum_ku_k(z)^*u_k(z)\,X(z)^{\frac{p}{2}}\|_1
 \le \|X(z)\|_p^p<1.$$
The same is true for $z\in\partial_1$. Hence, by the maximum
principle,
 $$\sum_k\|Y_k(\th)\|_2^2\le 1.$$
Set
 $$a=\Phi(\th),\quad y_k=Y_k(\th),\quad b=\Psi(\th).$$
Then by the previous discussion we have
 $$ay_kb=F_k(\th)=x_k\quad\mbox{and}\quad
 \|a\|_{r_0(\th)}\, \|b\|_{r_1(\th)}\,
 \big(\sum_k\|y_k\|_2^2\big)^{\frac{1}{2}}\le 1.$$
Thus $\a_{p,\th}(x)\le1$. This finishes the proof of Step 1.

\medskip

\n{\em Step~2: The proof of (\ref{interp1}) in the general case.}
Assume again $1\le p<2$; but now $M$ is a general von Neumann
algebra. The inequality $\|x\|_{p,\th}\le \a_{p,\th}(x)$ can be
proved as before by interpolation. Indeed, the same proof as in
Step~1 shows that this inequality still holds for $\th=0$ and
$\th=1$. Then we can use Terp's interpolation theorem
\cite{terp-int} to conclude as in Step~1. Alternatively, instead
of using Terp's theorem, we can also appeal to Kosaki's
interpolation theorem \cite{kos-int}, which is clearly applicable
to strictly normal semifinite faithful weights. $L_p(M)$ can be,
of course, constructed from such a weight on $M$.

\medskip

We will use Haagerup's reduction theorem to prove the converse
inequality. Keep all notations in Theorem \ref{haag}.  We consider
the noncommutative $L_p$-spaces based on $M$, $\M$ and $\M_i$.
$L_p(\M)$ is constructed with respect to the weight $\wh\f$ there.
Then i) and iv) of Theorem \ref{haag} imply that $L_p(M)$ and
$L_p(\M_i)$ can be considered, in a natural way, as (isometric)
subspaces of $L_p(\M)$. On the other hand, by ii), iv) and \cite[
Lemma~2.2]{jx-burk}, $\E$ and  $\E_i$ extend to contractive
projections from $L_p(\M)$ onto $L_p(M)$ and onto $L_p(\M_i)$,
respectively ($1\le p\le\8$). These extensions are still denoted
by the same symbols. Finally, by the w*-density of $\bigcup_i\M_i$
in $\M$ and \cite[Lemma~1.1]{jx-burk}, $\bigcup_i L_p(\M_i)$ is
dense in $L_p(\M)$ for $p<\8$. Moreover, by the commutation
relations $\E_i\circ\E_j=\E_j\circ\E_i$ in Theorem \ref{haag},
iv), we deduce that the family $\{L_p(\M_i)\}_{i\in I}$ is also
increasing, and for any $x\in L_p(\M)$ the net $\{\E_i(x)\}_i$
converges to $x$ in $L_p(\M)$ (relative to the w*-topology for
$p=\8$).

Note that $\E$ extends coordinate-wise to a projection from
$L_p(\M)^n$ onto $L_p(M)^n$, which is still denoted by $\E$. Then
$\E$ is contractive on $C_p^n[L_p(\M)]$ and on $R_p^n[L_p(\M)]$.
The same remark applies to each condition expectation $\E_i$ too.

By this complementation of $L_p(M)^n$ and $L_p(\M_i)^n$ in
$L_p(\M)^n$, we have the following isometric inclusions
 $$(C_p^n[L_p(M)],\; R_p^n[L_p(M)])_\th\subset
 (C_p^n[L_p(\M)],\; R_p^n[L_p(\M)])_\th$$
and
 $$(C_p^n[L_p(\M_i)],\; R_p^n[L_p(\M_i)])_\th\subset
 (C_p^n[L_p(\M)],\; R_p^n[L_p(\M)])_\th\;.$$
On the other hand, by the complementation of $L_p(\M_i)$  in
$L_p(\M_j)$ for $i\le j$ and the density of $\bigcup_i L_p(\M_i)$
in $L_p(\M)$,   we deduce  that the family
 $$\Big\{(C_p^n[L_p(\M_i)],\; R_p^n[L_p(\M_i)])_\th\Big\}_{i\in
 I}$$
is increasing and its union is dense in $(C_p^n[L_p(\M)],\;
R_p^n[L_p(\M)])_\th$. Consequently, for any $x\in L_p(\M)^n$,
$\{\E_i(x)\}_i$ converges to $x$ with respect to the interpolation
norm $\|\;\|_{p,\th}$.

Now fix $x\in L_p(M)^n$ with $\|x\|_{p,\th}<1$. Then  $\E_i(x)\in
L_p(\M_i)^n$, and by the previous discussion,
$\|\E_i(x)\|_{p,\th}<1$ for every $i\in I$. Since $\M_i$ is finite
and $\s$-finite,  we find, by Step~1, $a_i, b_i$ and
$(y_{k,i})_{1\le k\le n}$ in the unit balls of
$L_{r_0(\th)}(\M_i)$, $L_{r_1(\th)}(\M_i)$ and
$\el_2^n(L_2(\M_i))$, respectively, such that
 $$\E_i(x_k)=a_i\, y_{k,i}\, b_i\;, \quad 1\le k\le n,\ i\in I.$$
Define
 $$\el\;:\; L_{p'}(M)^n\to\comp\quad
 \mbox{by}\quad\el(z)=\sum_k\tr(x_k^*z_k).$$
Since $\E_i(x_k)\to x_k$ in $L_p(\M)$ for every $k$, we have
 $$\el(z)=\lim_i\sum_k\tr[\E_i(x_k)^*z_k]
 =\lim_i\sum_k\tr[y_{k,i}^*\,(a_i^*z_kb_i^*)],\quad
 \forall\; z\in L_{p'}(M)^n\;.$$
However, by Cauchy-Schwarz
 $$|\sum_k\tr[y_{k,i}^*\,(a_i^*z_kb_i^*)]\big|
 \le \big(\sum_k\|y_{k,i}\|^2_2\big)^{1/2}
 \big(\sum_k\|a_i^*z_kb_i^*\|^2_2\big)^{1/2}
 \le \a_{p',\th}(z).$$
It then follows that $\el$ is a contractive functional on
$(L_{p'}(M)^n,\; \a_{p',\th})$. It is clearly normal in the case
$p'=\8$. Therefore, by Lemma \ref{lemma} (and the remark following
it), we deduce $\a_{p,\th}(x)\le 1$ for $\el$ is defined by $x$.
This achieves the proof of (\ref{interp1}).

\medskip\n{\em Step~3: The proof of (\ref{interp2}).} Let $2<p\le\8$.
The inequality $\a_{p,\th}(x)\le\|x\|_{p,\th}$ can be easily
proved by interpolation. We omit the details. Let us prove the
converse by duality using (\ref{interp1}). Fix $x\in L_p(M)^n$.
Note that
 $$(C^n_{p'}[L_{p'}(M)],\; R^n_{p'}[L_{p'}(M)])_\th^*
 =(C^n_{p}[L_{p}(M)],\; R^n_{p}[L_{p}(M)])_\th \quad
 \mbox{isometrically}.$$
Thus by (\ref{interp1}) already proved
 \be
 \|x\|_{p,\th}
 &=&\sup\big\{\big|\sum_k\tr(y_k^*x_k)\big|\;:\;
 y\in L_{p'}(M)^n,\; \|y\|_{p',\th}\le 1\big\}\\
 &=&\sup\big\{\big|\sum_k\tr(y_k^*x_k)\big|\;:\;
 y\in L_{p'}(M)^n,\; \a_{p',\th}(y)\le 1\big\}
 \le\a_{p,\th}(x).
 \ee
This is the desired inequality. Therefore, the proof of Theorem
\ref{interp} is complete.\cqd


\end{document}